\documentclass[12pt]{amsart}

\usepackage{yhmath}

\usepackage{amssymb,amsmath,amsfonts,amscd,amsthm,enumerate}
 \usepackage{graphicx,color}
\graphicspath{{picturesfolder/}}
\usepackage{wrapfig}
\usepackage{multirow}
\usepackage{textcomp}

\usepackage{todonotes}

\usepackage{epsfig,subfigure}
\usepackage{todonotes}

\usepackage{lineno,xcolor}

\usepackage[all]{xy}
\usepackage{color}

\newcommand{\R}{{\Bbb R}}

\newcommand{\so}{\Rightarrow}
\newcommand{\e}{\varepsilon}

\newcommand{\inv}{^{-1}}

\newcommand{\del}{\partial}

\newtheorem{theorem}{Theorem}

\newtheorem{proposition}[equation]{Proposition}

\newtheorem{definition}[equation]{Definition}

\title{The Binary Returns!}

\author[Connor Jackman]{Connor Jackman}
\address[Jackman]{Mathematics Department, University of California,
4111 McHenry
Santa Cruz, CA 95064, USA}
\email{cfjackma@ucsc.edu}

\begin{document}

\maketitle
  
\begin{abstract}
Consider the spatial Newtonian three body problem at fixed negative energy and fixed angular momentum. The moment of inertia $I$ provides a measure of the overall size of a three-body system. We will prove that there is a positive number $I_0$ depending on the energy and angular momentum levels as well as the masses such that every solution  at these levels passes through $I\leq I_0$ at some instant of time. Motivation for this result comes from trying to prove the impossibility of realizing a certain syzygy sequence in the zero angular momentum problem.
\end{abstract}

\section{Introduction}
\label{intro}

The spatial 3-body problem concerns three point masses in space moving according to Newton's equations of gravitation. The point of this article is to prove that there exist no periodic solutions to this problem which ``hang out near infinity".

 The conserved quantities for the problem are the energy $H$, angular momentum $J$ and linear momentum. As is standard, we may, without loss of generality, assume that the linear momentum is zero and the origin of space coincides with the center of mass of the three bodies. If $m_i$ denote the masses and $q_i\in\R^3$ the positions of the bodies, then the standard measure of size is $\|q\|=\sqrt{I(q)}$ where $q=(q_1,q_2,q_3)$ and $I=\sum m_i |q_i|^2$ is known as the total moment of inertia. Neighborhoods of infinity are regions of the form $\{ q: I(q)\geq I_0\}$. As $I_0\to\infty$ the neighborhood converges to infinity. Our main theorem is:

\begin{theorem}  \label{mr} For $H<0$ there exists $I_0(m_j, H, J)>0$ such that any orbit at these energy and momentum levels beginning in the region $I>I_0$ enters the region $I\leq I_0$ in forwards or backwards time.

\end{theorem}

{\sc Motivation.}  The motivation behind  our result came from the problem of which syzygy sequences are realized in the \textit{zero angular momentum} planar three body problem (see \cite{homclass}, \cite{inf}, \cite{Syzygy}).
The term syzygy is from astronomy and refers to when the three bodies are in eclipse, that is collinear.  Each   syzygy has a `type' 1, 2, or 3, according to the  label of the mass in the middle.
Then the syzygy sequence of an orbit is this list of syzygy types  in temporal order.  A first open problem is
whether or not the periodic sequence of repeating 1212's is realized by a periodic solution to the zero angular momentum problem.
One imagines such a motion as consisting of masses 1 and 2 going around each other in a near circular orbit, very far from mass 3,
and the center of mass of the $m_1$ and $m_2$ orbit slowly going around mass 3, like the Earth-Moon-Sun system.
The action over such solutions decreases as the distance of the earth moon system to the sun goes to infinity i.e. minimizing the action forces the solution to slide off into a neighborhood of infinity, see \cite{ChM}.
The theorem excludes the existence of such solutions ``near infinity''  i.e in the region $I \ge I_0 (m_j, H, 0)$.\\

{\sc Remark.} In the theorem we may either exclude orbits having a binary collision singularity or we may pass through them using Levi-Civita regularization. Can we prove an analogous result to Theorem 1 for $N\ge 4$? The proof here breaks down in proposition \ref{claimapp} where the neighborhoods of infinity fail to split into connected components characterized by a far body with suitable Jacobi coordinates. The connectedness of the neighborhoods of infinity due to these spread out clusters of tight binaries is utilized for Jeff Xia's orbits realizing infinity in finite time singularities where $N\ge 5$. Can these infinity in finite time orbits provide counterexamples to Theorem 1 for $N\ge 5$?

{\sc Remark.} In \cite{Meyers} comet-like periodic orbits for the $N$-body problem are established in a region $I\ge I_C$ for $I_C$ large. These orbits do not contradict our theorem because as $I_C\to \infty$ their orbits angular momentum $|J|\to\infty$.

\section{Related Results.}
\label{sec0}

The behavior of $I(t)$ has long been studied to gain some qualitative understanding of the N-body problem. Sundman, \cite{Sundman}, showed for the three-body problem that non-zero angular momentum implies no orbits suffer triple collision i.e. $I>0$ for all orbits. Namely there exists a positive lower bound, $I_S(m_i, H, J, I(0), \dot I(0))$, for orbits at such levels. That is $I(t)>I_S>0$ over the solutions with energy $H$ and angular momentum $J\ne 0$ and with initial conditions at $I(0)$, $\dot I(0)$. Hadamard, \cite{Hadamard} pg. 259, gave an explicit formula for such an $I_S$ and G.D. Birkhoff \cite{Birkhoff} ch. IX \S8 studied escape conditions in the non-zero angular momentum case by showing for example (pg. 282) that $I$ sufficiently small (near zero)  at some instant, $t_0$, implies $I$ becomes infinite as $t$ goes to infinity. One might paraphrase Birkhoff's result as `no hanging out in neighborhoods of triple collision.'

A great deal of analysis on $I$ has followed these two tracks around small $I$ values. See for example \cite{LM}, \cite{MY} on the greatest lower bound of $I$ for bounded orbits and \cite{Marchal}, \cite{MY}, \cite{Escape} for  efficient tests of escape in a variety of cases. The book \cite{Marchal3b} ch. 11 is a detailed reference for the qualitative study of $I$.

For each orbit let $I_m$ be the minimum value of $I$ over this orbit, sharpening Sundman leads one to seek a (greatest) lower bound of $I_m$ over classes of orbits. An analogous question here is instead to seek a (least) upper bound of $I_m$ over all orbits.

While most focus in the literature so far appears on the greatest lower bound and escape this upper bound question has not entirely escaped notice. A statement similar to Theorem 1 appears in \cite{Marchal3b} pg. 468 where an upper bound is given in a remark about a class of equal mass cases (those with $H|J|^2=-\frac{4}{3^5}$) and the least upper bound is conjectured to be attained over the Brouke-Henon orbit (\cite{Marchal3b} pg. 469). Here we give a new motivation to this question as to the existence of the 1212... solution in the \textit{zero} angular momentum case and use a different method than that of \cite{Marchal3b}. Additionally we observe that both methods give upper bounds in a general case rather than just treating an equal mass case  (see also \cite{Marchal3b} pg. 483). Moreover the method we use here offers hope of lowering the upper bound if the perturbation step (propositions \ref{perclaim}, \ref{stripclaim}, \ref{Mainclaim}) is dealt with more effectively. In the appendix we give some comparison of the two methods.

\begin{figure}[ht!]
\centering
\includegraphics[width=90mm] {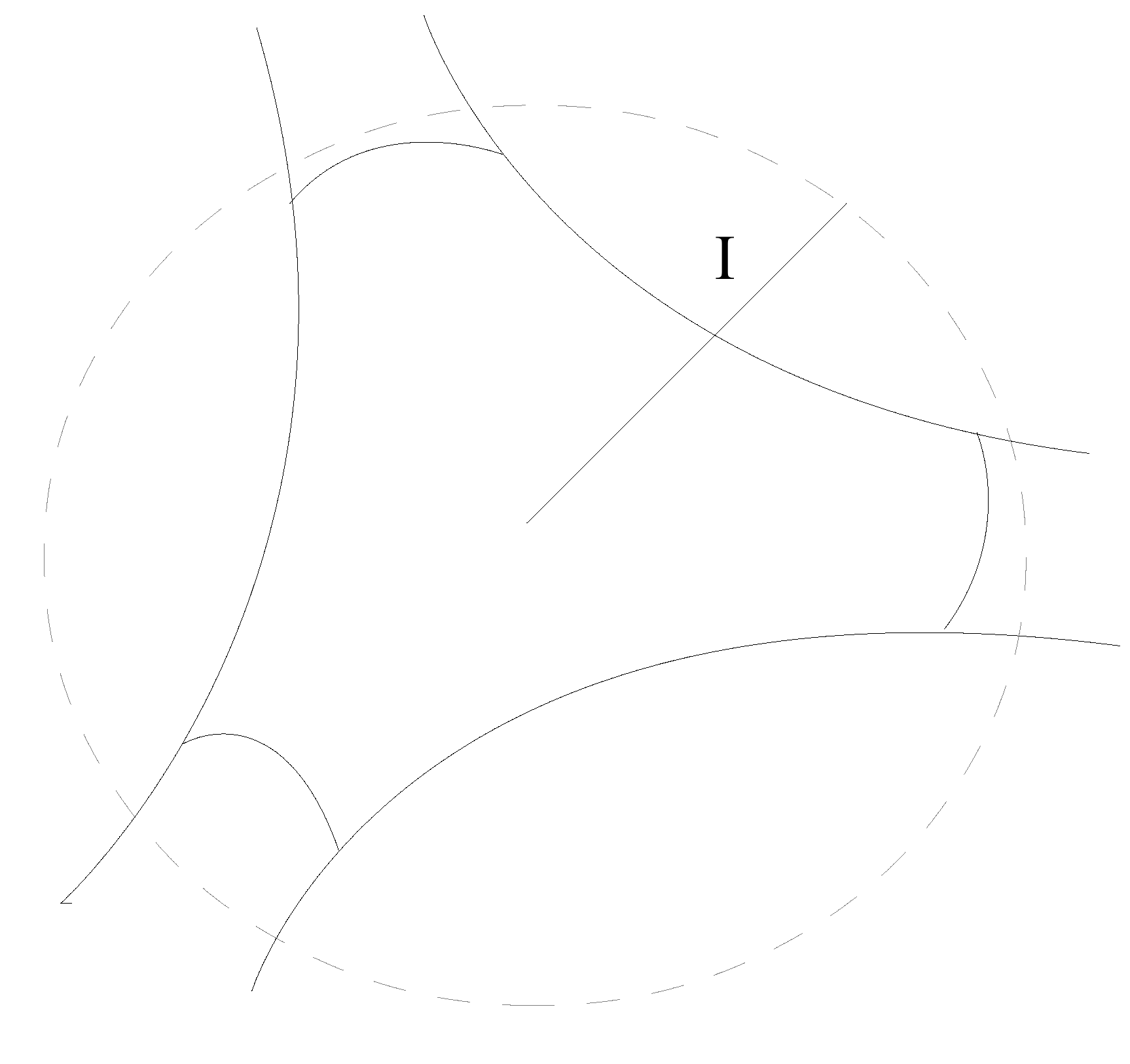}
\caption{For the planar 3-body problem the shape space is $\R^3$ where $I$ is the distance from the origin. The admissible configurations at fixed $H<0$ are interior to a pair of pants where each leg of the pants is asymptotic to a cylinder around a binary collision ray. See \cite{Moeck}, or \cite{Triangles} for details.}
\label{pants}
\end{figure}

\section{ Structure of proof.}
\label{sec1}

\begin{figure}[ht!]
\centering
\includegraphics[width=90mm] {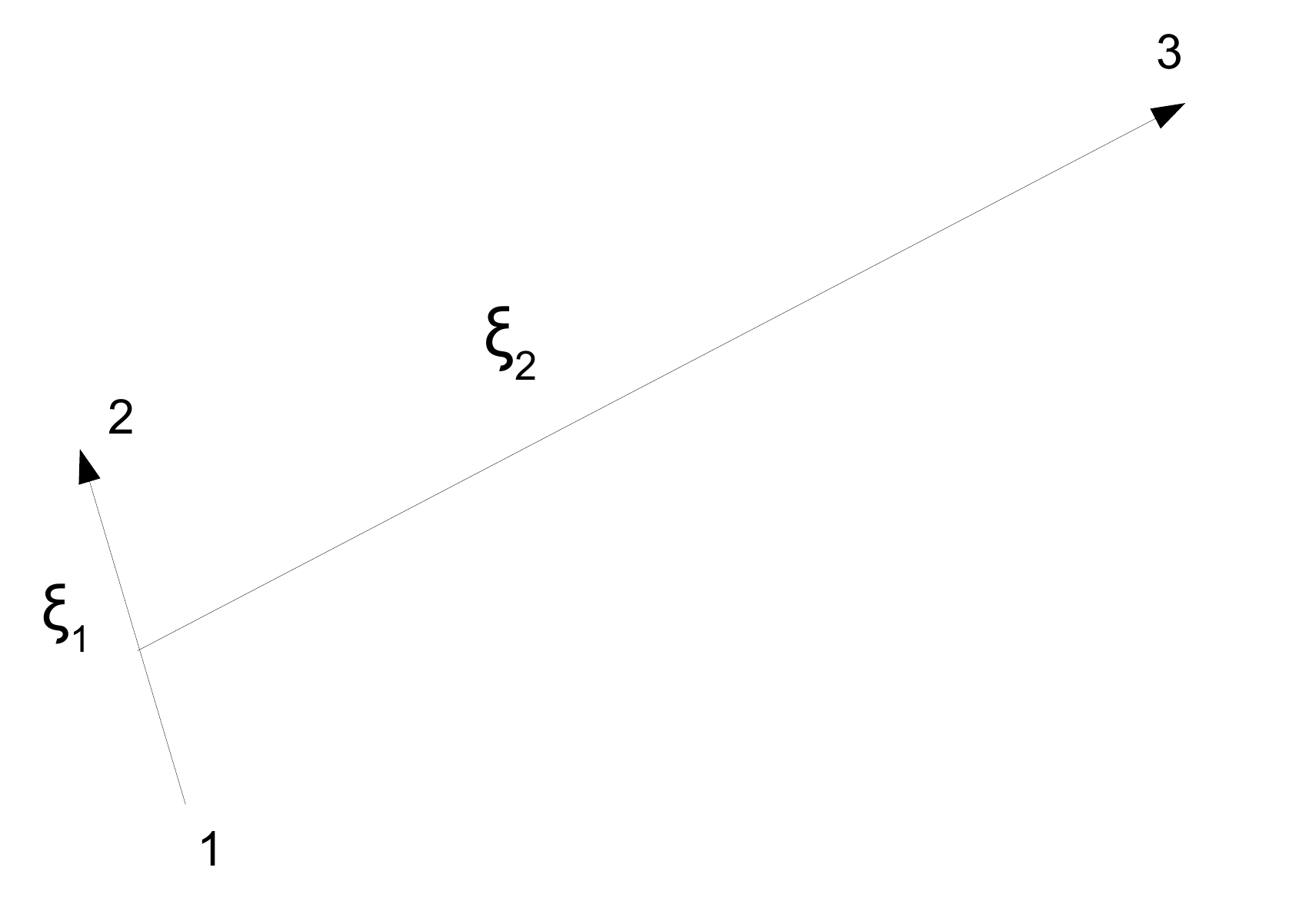}
\caption{A tight binary configuration, set $r:=|\xi_1|$, $\rho:=|\xi_2|$.}
\label{configs}
\end{figure}

For $H < 0$ as we let $I_0 $ increase, eventually the domain  $\{I \ge I_0\}$ splits into three components
 each  component   characterized by the  selection of one of the three masses.  The two remaining masses stay  close to each other while this   third selected  mass, stays relatively far away from
either member of this pair (see figure \ref{pants}).
We fix attention on one of these regions, supposing, after relabeling , that the close masses are $m_1$ and $m_2$.
In this region, we use the standard Jacobi coordinates $\xi_1, \xi_2$. See figure \ref{configs}.

When written in
these coordinates, Newton's differential equations becomes  a perturbation of two uncoupled Kepler problems,
one for each Jacobi vector, with the perturbation term getting arbitrarily small as $I_0 \to \infty$.   We focus attention
on the long Jacobi vector, which connects the center of mass of the $m_1$ and $m_2$ system to the 3rd mass.    When we
drop the perturbation term of this perturbed Kepler system, we get an exact solvable Kepler problem whose solutions we call
 ``the osculating solutions''.
 
 The Kepler parameters (energy, angular momentum, Laplace or Runge-Lenz vector) for the osculating system
can be {\it bounded } using that  $H, J$, the masses, are fixed
 and the fact that $I_0>>0$.   Now here comes the key observation, due to Chenciner.  Consider   a family of solutions to
 Kepler's equation having fixed energy and bounded angular momentum.  If,  along the solutions of this family
 the initial distance from the origin tends to infinity then these  orbits become extremely eccentric,  and thus
must come    close to the origin. Thus the osculating orbits   cannot ``hang out near infinity''.  Said  slightly differently,
since   large circular orbits
for the Kepler problem have large angular momentum and since our total angular momentum is fixed, large near circular motions
for osculating system are excluded and this excludes orbits of the type of our Earth-Moon-Sun cartoon described above.

Here is the strategy of proof then.  Show that for sufficiently large $I_0$ all of the osculating solutions starting in $\{ I \ge I_0\}$
are extremely eccentric, enough so to enter the region $\{ I\le I_0\}$ (see proposition \ref{oscclaim}).  Next show that the real solutions do not
vary too much from these osculating solutions,  as long as they stay in the region $I \ge I_0$, and for bounded times
(indeed for times of order $O(I_0 ^{3/2})$, proposition \ref{perclaim}). It   follows that if the osculating orbit enters the  region $I \le I_0$ within the time $O(I_0^{3/2})$ (which we expect by Kepler's third law) then the
true orbit must also enter into that region.
Finally, (proposition \ref{Mainclaim}) we verify that there is indeed sufficient time: the time scale over which the  approximation of the true motion
by the osculating motion is valid is long enough that   the true motions must follow their osculating leads into a region   $I\le I_0$.

\section{Set-up and Notation}
\label{sec2}

In the spatial 3-body problem, we consider the motion of three point masses $m_1,m_2,m_3$ under Newton's gravitational attraction. We will denote the configurations by $$q=(q_1,q_2,q_3)\in (\R^3)^3\backslash \{ (x_1,x_2,x_3): x_i=x_j\text{ some } i\neq j\}.$$

As is standard, we may take the center of mass zero coordinates ($\sum m_i q_i=0$) and will now define the Jacobi coordinates in which the splitting into two perturbed Kepler problems will be clear (see Figure \ref{configs} as well as \cite{Poll} 2.7, \cite{Fej}, or \cite{Stutt}):

$$\xi_1=q_2-q_1, $$ $$\xi_2=q_3-(m_1+m_2)\inv(m_1q_1+m_2q_2)=\frac{m_1+m_2+m_3}{m_1+m_2} q_3.$$

We set $$r=|\xi_1| \text{ and } \rho=|\xi_2|.$$

For reference we record here in one place the mass constants that will be used throughout:\\

{\sc Mass constants:}

$$\mu=m_1+m_2$$ $$M=m_1+m_2+m_3$$ $$\alpha_1=m_1m_2\mu\inv$$ $$\alpha_2=m_3\mu M\inv$$ $$\beta_1=\mu\alpha_1$$ $$\beta_2=M\alpha_2$$

Then in these coordinates we find:
\begin{equation}
I:=\sum m_i |q_i|^2=\alpha_1 r^2+\alpha_2\rho^2
\label{I}
\end{equation}
\begin{equation}
J:=\sum m_i (q_i\times\dot q_i)=\alpha_1 \xi_1\times\dot\xi_1+\alpha_2\xi_2\times\dot\xi_2=J_1+J_2
\label{J}
\end{equation}
for the moment of inertia, and angular momentum respectively. Also the energy splits into $$H=H_{kep}+g$$ where $$H_{kep}=\frac12\alpha_1 |\dot\xi_1|^2-\frac{\beta_1}{r}+\frac12\alpha_2 |\dot \xi_2|^2-\frac{\beta_2}{\rho}=H_1+H_2$$ is an energy for two uncoupled Kepler problems and $$g=\frac{\beta_2}{\rho}-\frac{m_1m_3}{|\xi_2+m_2\mu\inv\xi_1|}-\frac{m_2m_3}{|\xi_2-m_1\mu\inv\xi_1|}$$ is a perturbation term with $g=O(r^2/\rho^3)$, $g_{\xi_1}=O(r/\rho^3)$ and $g_{\xi_2}=O(r^2/\rho^4)$.

The equations of motion are then the two perturbed Kepler problems
\begin{equation}
\alpha_i\ddot \xi_i=-\frac{\beta_i\xi_i}{|\xi_i|^3}\mathbf{-}g_{\xi_i}
\label{peqs}
\end{equation}
\begin{definition}
A solution to the unperturbed Kepler problems satisfying the same initial conditions as a solution to these perturbed Kepler problems (eq. \ref{peqs}) will be called an  \textit{osculating orbit} (see \cite{Poll}, 1.16).
\end{definition}

\section{Proof of Main Theorem}
\label{sec3}

Fix the masses, angular momentum, negative energy $H<0$, linear momentum zero and a parameter $\lambda>0$ and only consider orbits at these energy and momentum levels in appropriate Jacobi coordinates. We will use $\overline I$ for a placeholder constant.


\begin{proposition}

For $H<0$, there exists $ I^*(m_i, H, J)>0$ such that the region $I> I^*$ consists of three connected components $B_1, B_2, B_3$. Moreover relabeling if necessary to fix our attention to $B_3$ (where $q_3$ is the far body) with appropriate Jacobi coordinates we have the bounds:
\begin{equation}
|g|\leq c_g(r^2/\rho^3),~~~|g_{\xi_2}|\leq c_{g_2}(r^2/\rho^4)
\label{perg}
\end{equation}
\begin{equation}
|J_2|\leq \alpha_2c_{J_2}
\label{Jbound}
\end{equation}
\begin{equation}
r\leq c_r
\label{rbound}
\end{equation}
on the perturbation term $g$ angular momentum $J_2$ and short Jacobi vector $r$ throughout $B_3$ for some constants $c_g, c_{g_2}, c_{J_2}, c_r$ depending on masses, energy and angular momentum.
\label{claimapp}
\end{proposition}

See \cite{Moeck}, \cite{Fej}, \cite{Stutt}, \cite{Marchal3b} regarding these well known \textit{lunar regions}.

\begin{proposition}

Take $I^{**}=\max\{ I^*, \alpha_1 c_r^2+\alpha_2  c_{J_2}^4/M^2\}$ where $I^*, c_r, c_{J_2}$ are from Proposition \ref{claimapp}. Then any osculating orbit with initial condition in $I> I^{**}$ falls in forwards or backwards time into the region $I\leq I^{**}$. Moreover the time to fall into the region $I\le I^{**}$ is less than or equal to the time to reach pericenter.

\label{oscclaim}
\end{proposition}

Proof: By eqs. (\ref{I}, \ref{rbound}) in the region $I>I^{**}$ we have $\rho^2> c_{J_2}^4/M^2$.

The `$\rho$' component of the osculating orbit of an initial condition in $I> I^{**}$ is a solution to the Kepler problem $$\ddot\xi_{osc}=-M\xi_{osc}/|\xi_{osc}|^3$$ with $\rho_{osc}^2(0)=|\xi_{osc}(0)|^2>c_{J_2}^4/M^2$ and the restriction from eq. (\ref{Jbound}) $$|\xi_{osc}\times \dot\xi_{osc}|=\alpha_2\inv |J_2(0)|\leq c_{J_2}$$ on the angular momentum. Also from Proposition \ref{claimapp}, we have the $r$ component satisfying $r\leq c_r$ as long as we remain in the region $I>I^*$.

We now verify that for all such orbits, $\xi_{osc}$, the pericenter distance, $\rho_{osc}^{pc}$ is bounded.\\

case 1: $J_2\neq 0$.

In polar coordinates, any non-collision osculating orbit is (for some $e\geq 0$):
$$\rho_{osc}=\frac{\alpha_2^{-2}|J_2(0)|^2}{M(1+e\cos\theta)}$$ where $\theta=0$ corresponds to the pericenter.

Then as $e\geq 0$ and by eq. (\ref{Jbound}), $$\rho_{osc}^{pc}=\frac{\alpha_2^{-2}|J_2(0)|^2}{M(1+e)}\leq \frac{c_{J_2}^2}{M}.$$

case 2: $J_2=0$.

Collision! So the pericenter distance in this case is zero.\\

Now an osculating orbit starting in $I> I^{**}$ either reaches pericenter or leaves $I^*$ before it reaches pericenter. If it reaches pericenter before leaving $I>I^*$ then we have $I_{pc}\leq \alpha_1c_r^2+\alpha_2 c_{J_2}^4/M^2\le I^{**}$ so in either case we fall into the region $I\leq \max\{I^*, \alpha_1 c_r^2+\alpha_2 c_{J_2}^4/M^2\}= I^{**}$ in forwards or backwards time which is no more than $t_{pc}$, the time to pericenter. \qed

\begin{proposition}

Let $\overline I\ge \max\{I^*, \alpha_1c_r^2+\max\{1,(\frac{3c_{J_2}^2}{2M})^2\}\alpha_2 \}=\overline R$. Set $\overline\rho=\sqrt{\alpha_2\inv (\overline I-\alpha_1 c_r^2)}$ and $\e=1/\overline\rho$. Then any orbit with initial condition in $I\ge \overline I$ satisfies:
\begin{equation}
|\rho(t)-\rho_{osc}(t)|<A_1\e
\label{rhovar}
\end{equation}
for time
\begin{equation}
|t|\leq B_1\e^{-3/2}
\label{time}
\end{equation}
throughout the region $I\geq \overline I$.

Here we may pick the constant $B_1>0$ and then define $A_1=\frac{a}{M}(2+e^{\sqrt{2M+3c_{J_2}^2}B_1})$ where $a=\alpha_2\inv((c_{g_2}c_r^2 B_1)^2+2c_{J_2}c_{g_2}c_r^2 B_1+c_{g_2}c_r^2)$.
\label{perclaim}
\end{proposition}

Proof: First, from eq. (\ref{rbound}) any configuration with  $I\ge\overline I$ has $\mathbf{\rho}\ge\overline\rho\ge \max\{ 1, \frac{3c_{J_2}^2}{2M}\}\ge \max\{ 1, \frac{3\alpha_2^{-2}|J_2|^2}{2M}\}$, in particular our initial condition.

We consider our perturbed Kepler problem for the `$\rho$' motion:

$$\ddot\xi_2=-\frac{M\xi_2}{\rho^3}+F(\xi_2, t)$$

Where the time dependence in the perturbation term $F=-\alpha_2\inv g_{\xi_2}$ is due to the interaction of the motion of masses 1 and 2.

In the region $I\ge\overline I$ , we have $|F|\leq \alpha_2\inv c_{g_2}c_{r}^2\rho^{-4}\leq \alpha_2\inv c_{g_2}c_{r}^2\e^4$. We will set $$A=\alpha_2\inv c_{g_2}c_{r}^2.$$

An estimate for the variation of $c_t^2:=|\xi_2\times\dot\xi_2|^2=\alpha_2^{-2} |J_2(t)|^2$ will be needed. Since $|\dot c|\leq |\alpha_2\inv\dot J_2|=|\xi_2\times F|\leq A\rho^{-3}$, we have $$|\dot c|\leq A\e^3$$ so that
$$|c_t-c_0|\leq A\e^3 |t|.$$

Hence $$|c_t^2-c_0^2|\leq A\e^3|t|(A\e^3|t|+2c_0)\leq A\e^3|t|(A\e^3|t|+2c_{J_2}).$$

so that for $|t|\leq B_1\e^{-3/2}$ and $I\ge\overline I$ with $b=(AB_1)^2+2c_{J_2}AB_1$ we have
\begin{equation}
|c_t^2-c_0^2|\leq b\e^{3/2}.
\label{c2var}
\end{equation}

provided $\e\le 1$ which is guaranteed so long as $\overline I\ge \alpha_1 c_r^2+\alpha_2$ as is indeed the case since $\overline I\ge \overline R$.

To prove the proposition we'll use the Sandwich Lemma (see \cite{Syzygy} pg. 1942. Note that in \cite{Syzygy} there is an unneeded assumption requiring that $F_+<0$):\\

Sandwich Lemma: Given $\ddot x_{-}=F_{-}(x_{-})$, $\ddot x=F(x,t)$ and $\ddot x_{+}=F_{+}(x_{+})$ satisfying $F_{-}(x)\leq F(x,t)\leq F_{+}(x)$ and $\frac{\del F_{\pm}}{\del x_{\pm}}\geq 0$ over some time interval, then over this same time interval the solutions to $F_{\pm}, F$ satisfying the same initial conditions have:

$$x_{-}(t)\leq x(t)\leq x_{+}(t).$$\\

Now: $$\rho_{osc}\ddot\rho_{osc}+\dot\rho_{osc}^2=\frac{d}{dt} \rho_{osc}\dot\rho_{osc}=\frac{d}{dt} \xi_{osc}\cdot\dot\xi_{osc}=-M\rho_{osc}^2\rho_{osc}^{-3}+|\dot \xi_{osc}|^2 =-M\rho_{osc}\inv+\dot\rho_{osc}^2+c_0^2\rho_{osc}^{-2}$$ so

$$\ddot\rho_{osc}=c_0^2\rho_{osc}^{-3}-M\rho_{osc}^{-2}.$$

And likewise:

\begin{equation}
\ddot\rho= c_t^2\rho^{-3}-M\rho^{-2}+f(t)
\label{2der}
\end{equation}

where $|f(t)|=|\rho(t)\inv (\xi_2(t)\cdot F(\xi_2(t), t))|\leq A\rho(t)^{-4}$.\\

Take $v_1(\rho)=c_0^2 \rho^{-3}-M\rho^{-2}$ and $v_2(\rho,t)=c_t^2 \rho^{-3}-M\rho^{-2}+f$. We view $f$ here as $f(t)$ by plugging the true solutions $\xi_1(t), \xi_2(t)$ into $F, \rho$.

Now using our $|c_t^2-c_0^2|$ estimate eq. (\ref{c2var}) and our bound on $f$ we get:
$$|v_1-v_2|\leq b\e^{9/2}+A\e^4\leq a\e^4$$ for $a=b+A$, or
$$v_1-a\e^4\leq v_2\leq v_1+a\e^4$$

for time $|t|\leq B_1\e^{-3/2}$ and $I\ge \overline I$.\\

Now $\rho$ is a solution to $\ddot \rho=v_2$ and let $\rho_{\pm}$   be solutions to: $$\ddot \rho_\pm=v_1(\rho_{\pm})\pm a\e^4=:F_{\pm}(\rho_\pm)$$ satisfying the same initial conditions as $\rho$. Throughout the region $I\ge\overline I$ we have $\rho_{\pm}\ge\frac{3c_0^2}{2M}$ which implies  $\frac{\del F_{\pm}}{\del \rho_{\pm}}\ge 0$ so that we may apply the Sandwich Lemma throughout the region $I\ge\overline I$ yielding:

$$\rho_{-}\leq \rho\leq \rho_{+}$$ for time $|t|\leq B_1\e^{-3/2}$ as long as we remain in the region $I\ge \overline I$.\\

Likewise since $v_1-a\e^4\leq v_1\leq v_1+a\e^4$, we have for $|t|\leq B_1\e^{-3/2}$ and throughout $I\ge\overline I$ that
$$\rho_{-}\leq \rho_{osc}\leq \rho_{+}$$ holds.

Now we will show $\rho_+$ and $\rho_{-}$ remain close to finish the proof. Set $\eta=\rho_+-\rho_{-}\geq 0$.

Note that $v_1$ is Lipschitz in the region $\rho\ge \overline\rho$ with

$|v_1(x)-v_1(y)|\leq \omega|x-y|$ for $x,y\ge\overline\rho$ and $\omega=(2M+3c_0^2)\e^3=k\e^3$.\\

Then $\ddot \eta=v_1(\rho_+)-v_1(\rho_{-})+2a\e^4\so |\ddot \eta|\leq \omega |\eta|+2a\e^4=\omega\eta +2a\e^4$ so $$|\ddot\eta|\leq \omega \eta +2a\e^4.$$

Let  $F=v_1(\rho_+)-v_1(\rho_{-})+2a\e^4$ then we have $0\le F\le \omega\eta+2a\e^4$ provided $\rho_{-}\le \rho_+$ and $\frac{\del v_1}{\del\rho}>0$; which indeed holds throughout the region $I\ge\overline I$ for time $|t|\le B_1\e^{-3/2}$. Now the Sandwich Lemma with $F_{+}(\eta)=\omega\eta+2a\e^4$ and $F_{-}=0$ gives:

$$0\leq \eta(t)\leq \frac{2a\e^4}{\omega}(\cosh\sqrt\omega t-1)$$

and since $\omega=k\e^3$ where $2M\le k\le 2M+3c_{J_2}^2$ we have

$$|\rho(t)-\rho_{osc}(t)|\leq \rho_+(t)-\rho_{-}(t)=\eta(t)\leq \frac{2a\e}{k}(2+e^{\sqrt\omega |t|})\leq A_1\e$$

for time $|t|\leq B_1\e^{-3/2}$ as long as we are in the region $I\ge\overline I$ and where we set $A_1=\frac{a}{M}(2+e^{\sqrt{2M+3c_{J_2}^2} B_1})$.\qed

\begin{proposition}
Set $R=\max\{ \overline R, I^{**}, 4\alpha_1c_r^2\}$  where $\overline R$ is from proposition \ref{perclaim}. For $\overline I\ge R$ set $$\overline I^+=4(\overline I-\alpha_1c_r^2)>\overline I.$$ Then for any orbit with an initial condition in the strip $$\overline I\leq I\leq \overline I^+$$  we have that eq. (\ref{rhovar}) holds with $B_1=2^{3/2}\pi \sqrt{M}$ until the osculating orbit enters the region $I\le \overline I$.
\label{stripclaim}
\end{proposition}

Proof: First consider orbits with initial condition in $I\ge \overline I$ for some $\overline I\geq \max\{I^{**}, \overline R\}$ and with $ \e=1/\overline \rho$ defined as in proposition \ref{perclaim} and recall that $I\geq \overline I$ implies that $\rho\geq\overline\rho$. For osculating collision orbits with $J_2(0)=0$, some energy $H_2$ and $\rho(0)=\rho_{osc}(0)>\overline\rho$ the time to collision in forwards time (or time from expulsion in backwards time) $t_c$, satisfies:
\begin{equation}
t_c\leq \pi (8M)^{-1/2}\rho_{osc}(0)^{3/2}.
\label{coll}
\end{equation}

We will use Lambert's Theorem (see \cite{rec}) to compare time to pericenter for general osculating orbits to these collision times. Lambert says that for Kepler orbits, the time of travel between two points, $a_1, a_2$ on the orbit is a function of the energy, chord length $d=|a_1-a_2|$ and $|a_1|+|a_2|=r_1+r_2$ (where the origin is at the focus, see figure \ref{lambert}). Namely, for equivalent configurations (those having the same energy, same chord length $d$, and $r_1+r_2=s_1+s_2$) the time of travel from $a_2$ to $a_1$ is the same as the time of travel from $b_2$ to $b_1$. Figure \ref{lambert} is how we will choose our equivalent configurations:

\begin{figure}[ht!]
\centering
\includegraphics[width=90mm] {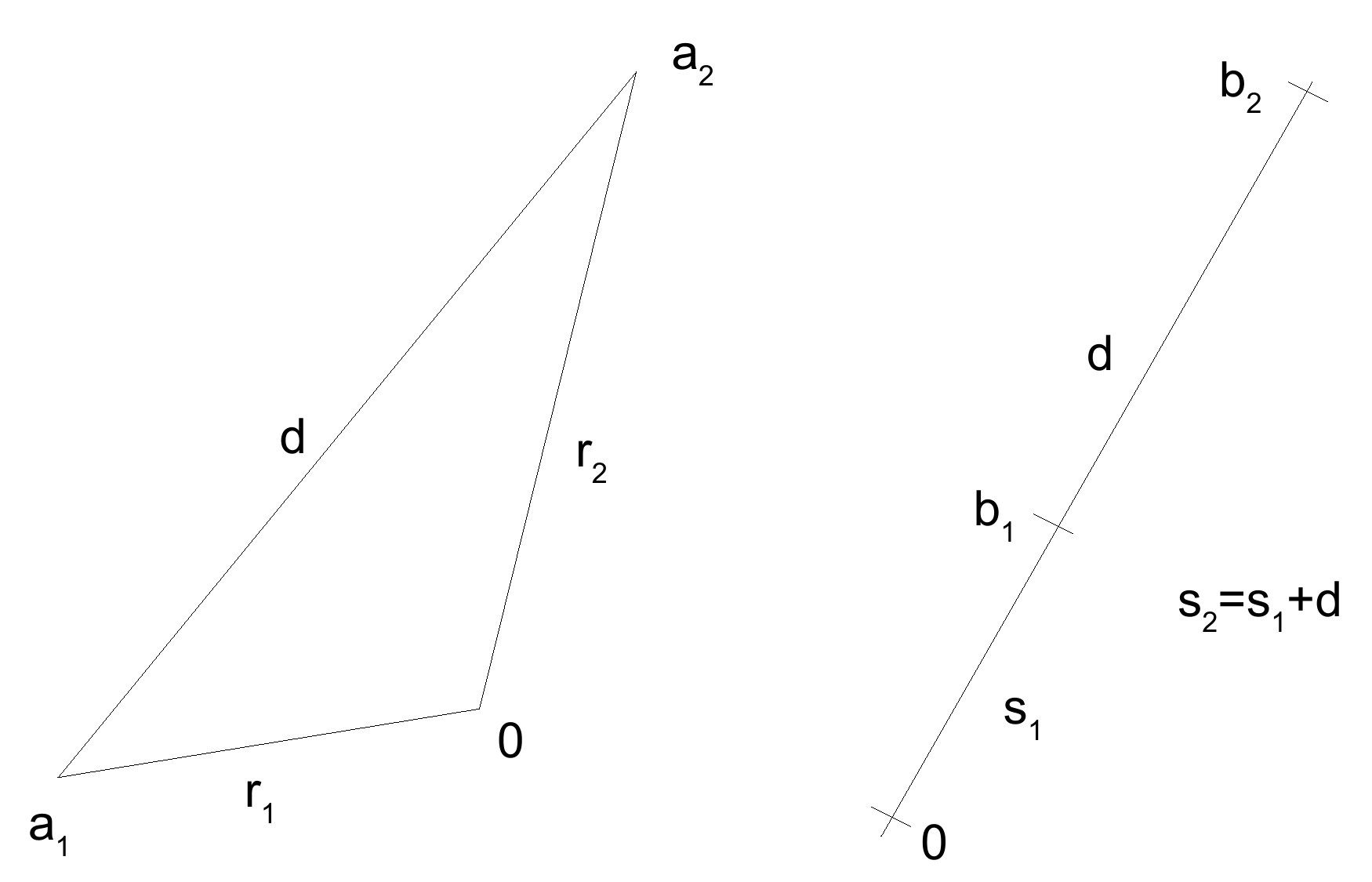}
\caption{Two equivalent configurations.}
\label{lambert}
\end{figure}

For a general osculating orbit $\rho_{osc}$, take $r_1=\rho_{osc}^{pc}$, $r_2=\rho_{osc}(0)=\rho(0)>\overline\rho$ and then $s_1, s_2$ are determined by $s_2-s_1=d=|a_2-a_1|$ and $s_1+s_2=r_1+r_2$. By Lambert's theorem and eq. (\ref{coll}) we have the time to pericenter, $t_{pc}$, satisfies
\begin{equation}
t_{pc}\leq \pi(8M)^{-1/2} s_2^{3/2}
\label{l1}
\end{equation}
And since $r_2\ge r_1$ (as we are in $I\ge I^{**}$) we have:
$$2s_2-(r_1+r_2)=s_2-s_1=d\le r_1+r_2\so $$ $$s_2\leq r_1+r_2\le 2r_2.$$
So continuing with eq. (\ref{l1}) $$t_{pc}\leq \pi M^{-1/2}r_2^{3/2}.$$
So to compare  $t_{pc}$ with our estimates eq. (\ref{time}) we want  $t_{pc}\leq B_1\e^{-3/2}=B_1\overline\rho^{3/2}$, which holds when:
$$\pi M^{-1/2}r_2^{3/2}\leq B_1\overline\rho^{3/2}\so $$
$$r_2^{3/2}\leq \pi\inv M^{1/2}B_1\overline\rho^{3/2}.$$

Take $B_1=2^{3/2}\pi/\sqrt{M}$ so that we will be working in the strip:

$$\overline\rho^{3/2}\leq r_2^{3/2}\leq 2^{3/2}\overline\rho^{3/2}$$
i.e. (recall that $r_2=\rho(0)=\rho_{osc}(0)$)

$$\overline\rho\leq \rho\leq 2\overline\rho.$$

The condition $\rho\leq 2\overline\rho$ is ensured (eqs. \ref{I}, \ref{rbound}) when  $I\le\overline I^+:=4\alpha_2\overline\rho^2=4(\overline I-\alpha_1 c_r^2)$.

Also, we ensure $\overline I<\overline I^+=4(\overline I-\alpha_1c_r^2)$ provided $\overline I\ge 4 \alpha_1 c_r^2 > \frac43 \alpha_1 c_r^2$.\qed

\begin{proposition}
(Main Theorem) Fix a parameter $\lambda>0$. Then there exists $R_\lambda(m_i, H, J)>0$ such that any orbit with initial condition satisfying $I(0)\ge R_\lambda$ comes in forwards or backwards time into the region $I\le R_\lambda$.

Explicitly, take $\overline R_\lambda=\max\{R,  \alpha_1c_r^2+\alpha_2(\frac{\alpha_2 A_1^2}{\lambda}), 2\alpha_2 A_1+4\alpha_1 c_r^2 +\lambda\}$ and $R_\lambda=\overline R_{\lambda}+2\alpha_2A_1+\lambda$   where $R$ is from proposition \ref{stripclaim}.
\label{Mainclaim}
\end{proposition}

Proof: Take $\overline I\ge \overline R_\lambda$ and $\e\inv=\overline\rho=\sqrt{\alpha_2\inv(\overline I-\alpha_1c_r^2)}$ and consider an orbit with initial condition in $\overline I^+ \ge I\geq\overline I$ as in proposition \ref{stripclaim}. By proposition \ref{oscclaim} we can let $t^*$ be the time the osculating orbit hits $\alpha_1 c_r^2+\alpha_2\rho_{osc}^2( t^*)=\overline I$ i.e. $\rho_{osc}(t^*)=\overline \rho=\e\inv$.\\

Along the true motion then at $t^*$ we have by proposition \ref{stripclaim} (eqs. \ref{rbound}, \ref{rhovar}) that

$$I(t^*)\le \alpha_1 c_r^2+\alpha_2\rho(t^*)^2\leq \alpha_1c_r^2+\alpha_2(\rho_{osc}(t^*)+A_1\e)^2=\overline I+2\alpha_2 A_1+\alpha_2 A_1^2\e^2$$

holds. Moreover due to the condition $\overline I\ge \overline R_\lambda\ge \alpha_1 c_r^2+\alpha_2(\frac{\alpha_2A_1^2}{\lambda})$ we have $\e^2\le \frac{\lambda}{\alpha_2 A_1^2}$ so that $$I(t^*)\le \overline I+2\alpha_2A_1+\lambda.$$ Also the condition $\overline I\ge \overline R_\lambda\ge 2\alpha_2 A_1+4\alpha_1 c_r^2 +\lambda>\frac13(2\alpha_2 A_1+4\alpha_1 c_r^2 +\lambda)$ ensures that $\overline I+2\alpha_2A_1+\lambda <\overline I^+.$

That is taking any $\overline I\ge \overline R_\lambda$ and setting $\lambda'=2\alpha_2A_1+\lambda$ then all orbits with initial condition in the strip $$\overline I+\lambda'\leq I(0)\le \overline I^+$$ come in forwards or backwards time into the region $$I\le \overline I+\lambda'$$.

\begin{figure}[ht!]
\centering
\includegraphics[width=90mm] {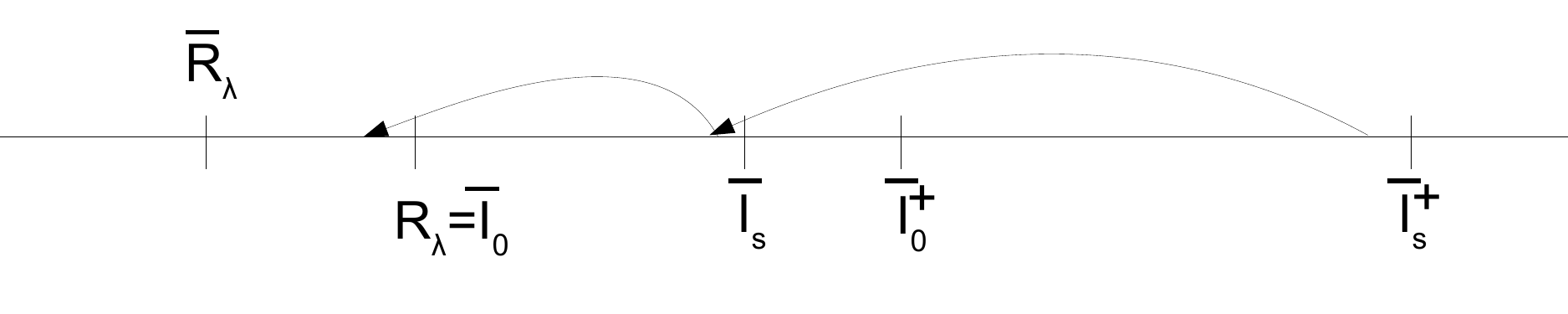}
\caption{Jumping back along the strips.}
\label{jump}
\end{figure}

In particular by setting $$\overline I_s=\overline R_\lambda+s$$ for $s\ge 0$, we may exhaust the region $I\ge R_\lambda=\overline R_\lambda+\lambda'=\overline I_0+\lambda'$ with the strips $$\overline I_s+\lambda'\le I\le \overline I_s^+.$$ Note that $s>s'$ implies $\overline I_s^+-\overline I_s>\overline I_{s'}^+-\overline I_{s'}$. Hence any orbit with initial condition $I(0)\ge R_\lambda$ will be forced to jump back along the strips (see figure \ref{jump}).

Finally in Theorem \ref{mr} we can take $I_0=R_{\lambda}$ for any choice of $\lambda$ (for instance $I_0=\min_{\lambda\in (0,1)} R_\lambda$). \qed

\section{Acknowledgments}
I would like to thank Richard Montgomery for many patient discussions, guidance with proofreading and writing, and e-mail introductions to Alain Chenciner, who had the original idea and question, as well as Ken Meyer and Rick Moeckel who gave helpful references. Alain Albouy also gave valuable help with the references to related results.

\begin{appendix}
\label{app2}

\section{Comparison to Marchal's equal mass case}

We shall now compare our results to \cite{Marchal3b} by examining the case: $$m_i=\frac13,~H=-\frac16,~ |J|=\frac{\sqrt{8}}{9}.$$

As $I^*$ corresponds to the apocenter of the collinear Euler motion (where $U|_{I=1}$ has a saddle point), we have $$I^*=\frac{32}{27}.$$

Here Marchal observed in \cite{Marchal3b} pg. 468 that $\ddot \rho<0$ for $I>I_M$ where $$I_M=\frac13(2.709629...)^2=2.447363...$$ so that every orbit will enter the region $I\le I_M$ at some instant (of course excluding or passing through any binary collisions). Moreover it was conjectured that in fact all orbits pass below the minimal inertia of the Henon-Brouke orbit (see \cite{Marchal3b} pg. 469) which is approximately 2.402035....and resembles an $e=0$ earth-moon sun cartoon type orbit.

Applying our final result (proposition \ref{Mainclaim}) in this case we obtain an $I_0>I_M$. However our bound only becomes larger when we apply our perturbation arguments; in this case and in general we obtain a lower pre-perturbation $I^{**}<I_M$ (of proposition \ref{oscclaim}). This gives hope that if our perturbation methods are improved (propositions \ref{perclaim}, \ref{stripclaim}, \ref{Mainclaim}) then the Marchal's bound of $I_M$ could be lowered. Now we outline how $I^{**}<I_M$ in general.

Note that Marchal's observation on the negativity of $\ddot\rho$ works not just in this equal mass case but lends to a shorter proof of Theorem 1 by using eq. \ref{2der}. We set $\mu_i=m_i\mu^{-1}$, $\lambda=r/\rho$, $\gamma=\angle(\xi_1, \xi_2)$ and rewrite eq. \ref{2der} as $$\rho^3\ddot\rho=c^2-\rho M\phi(\lambda,\gamma)$$ where $$\phi(\lambda,\gamma)=\mu_1\frac{1+\mu_2\cos(\gamma)\lambda}{(1+2\mu_2\cos(\gamma)\lambda+\mu_2^2\lambda^2)^{3/2}}+\mu_2\frac{1-\mu_1\cos(\gamma)\lambda}{(1-2\mu_1\cos(\gamma)\lambda+\mu_1^2\lambda^2)^{3/2}}.$$

Note that $\phi\ge\delta>0$ throughout $I\ge I^*$ for a $\delta\le \phi(\lambda, \pi/2)<1$ dependent on the masses. Hence a $\rho_M$ corresponding to Marchal's $I_M$ is (in general) $$\rho_M=\frac{c_{J_2}^2}{M\delta}.$$

However although eq. \ref{2der} leads to a simpler proof, following an orbit to pericenter rather than over the region where $\ddot\rho<0$ has the potential to yield lower upper bounds as for Keplerian orbits we have $$\rho^{pc}=\frac{c^2}{M(1+e)}\le \frac{c_{J_2}^2}{M}<\rho_M.$$ Thus the $I^{**}$ of proposition \ref{oscclaim} satisfies $$I^{**}<I_M.$$

So our strategy of proof provides hope of lowering the bound towards Marchal's  conjectured Henon-Broucke value in this case and a lower upper bound in general \textit{provided} the techniques in the perturbation steps propositions \ref{perclaim}, \ref{stripclaim}, \ref{Mainclaim} are improved to follow the orbits past the $\ddot \rho<0$ regime. Can they be improved? Perhaps in some non-equal mass cases or for some (outer) eccentricity $e>0$ orbits above Henon-Broucke? I am optimistic that taking advantage of the sharper bounds and techniques of the literature they can be improved at least for large classes of orbits. Especially so as many bounds of the perturbation steps here are not the sharpest (as the original motivation here was merely the existence of some upper bound in general specifically the zero angular momentum case).

\end{appendix}




\end{document}